\documentclass[journal]{IEEEtran}

\usepackage{cite}

\ifCLASSINFOpdf
  \usepackage[pdftex]{graphicx}
\else
\fi
\usepackage{amsmath}
\usepackage{amssymb}
\usepackage[short]{optidef}
\usepackage{bm}
\makeglossary
\DeclareMathOperator{\E}{E}
\DeclareMathOperator{\sign}{sign}
\usepackage{siunitx}

\usepackage[utf8]{inputenc}
\usepackage[gen]{eurosym}
\DeclareUnicodeCharacter{20AC}{\euro{}}

\usepackage{pgfplots,tikz}
\pgfplotsset{compat=newest} 
\pgfplotsset{plot coordinates/math parser=false} 
\newlength\figureheight 
\newlength\figurewidth 
\usepgfplotslibrary{external} 
\usepackage{url}

\usepackage[hidelinks]{hyperref}
\makeatletter
\g@addto@macro\UrlBreaks
{%
	\do\a\do\b\do\c\do\d\do\e\do\f\do\g\do\h\do\i\do\j%
	\do\k\do\l\do\m\do\n\do\o\do\p\do\q\do\r\do\s\do\t%
	\do\u\do\v\do\w\do\x\do\y\do\z\do\&\do\1\do\2\do\3%
	\do\4\do\5\do\6\do\7\do\8\do\9\do\0\do\/\do\.%
}
\g@addto@macro\UrlSpecials
{%
	\do\/{\mbox{\UrlFont/}\hskip 0pt plus 10pt}%
}

\usepackage[shortcuts]{extdash}

\begin{document}
%
\title{Combined Stochastic Optimization of Frequency Control and Self-Consumption with a Battery}
%
%
%

\author{Jonas~Engels, 
        Bert~Claessens 
        and~Geert~Deconinck

\thanks{Jonas Engels is with REstore, Antwerp, Belgium and with the Department of Electrical Engineering, KU Leuven/EnergyVille, Leuven, Belgium  (jonas.engels@restore.energy)}
\thanks{Bert Claessens is with REstore, Antwerp, Belgium (bert.claessens@restore.energy)}
\thanks{Geert Deconinck is with the Department of Electrical Engineering, KU Leuven/EnergyVille,	Leuven, Belgium (geert.deconinck@kuleuven.be)}
\thanks{This work is partially supported by Flanders Innovation \& Entrepreneurship (VLAIO)}}

%


\maketitle

%

\begin{abstract} Optimally combining frequency control with self-consumption can increase revenues from battery storage systems installed behind-the-meter. This work presents an optimized control strategy that allows a battery to be used simultaneously for self-consumption and primary frequency control. Therein, it addresses two stochastic problems: the delivery of primary frequency control with a battery and the use of the battery for self-consumption.

We propose a linear recharging policy to regulate the state of charge of the battery while providing primary frequency control. 
Formulating this as a chance-constrained problem, we can ensure that the risk of battery constraint violations stays below a predefined probability. We use robust optimization as a safe approximation to the chance-constraints, which allows to make the risk of constraint violation arbitrarily low, while keeping the problem tractable and offering maximum reserve capacity. Simulations with real frequency measurements prove the effectiveness of the designed recharging strategy.

We adopt a rule-based policy for self-consumption, which is optimized using stochastic programming. This policy allows to reserve more energy and power of the battery on moments when expected consumption or production is higher, while using other moments for recharging from primary frequency control. We show that optimally combining the two services increases value from batteries significantly.
\end{abstract}


%
\IEEEpeerreviewmaketitle

\section*{Nomenclature}
\addcontentsline{toc}{section}{Nomenclature}
\begin{IEEEdescription}[\IEEEusemathlabelsep\IEEEsetlabelwidth{$c_r, c_{cons}, c_{inj}$}]
	\item[$c_r, c_{cons}, c_{inj}$] Price of frequency control capacity, consumed and injected energy
	\item[CVaR] Conditional value-at-risk
	\item[$D$] Recharge policy matrix
	\item[$E$] Energy content
	\item[$\E {[\cdot ]}$] Expected value
	\item[$f$] Frequency
	\item[$g$] Objective function
	\item[$i, j, k$] Index, scenario index, time step
	\item[$I$] Identity matrix
	\item[$n_c, n_t, n_{sc}$] Number of constraints, number of time steps, number of scenarios
	\item[$P$] Power
	\item[$\Pr$] Probability
	\item[$r$] Frequency control capacity
	\item[$SoC$] State of Charge
	\item[$W$] Whitening transformation matrix
	\item[$\alpha$] Confidence bound
	\item[$\Delta f$] Frequency deviation from 50 Hz
	\item[$\Delta t$] Time step duration
	\item[$\epsilon$] Probability of constraint violation
	\item[$\lambda$] Dual variables
	\item[$\eta$] Efficiency
	\item[$\sigma_f, \sigma_b$] Forward/backward deviations
	\item[$\Sigma_{\bm{x}}$] Covariance matrix of $\bm{x}$
	\item[$\mathcal{F}$] Uncertainty set 
	\item[$\overline{[\cdot]}$] Mean value
	\item[${[\cdot]}^c, {[\cdot]}^d$] Charging, discharging
	\item[${[\cdot]}_k, {[\cdot]}_i$] Value at time step $k$, index $i$
	\item[$\bm{x}$] A bold symbol indicates the vector of the values $(x_1, x_2, \ldots, x_{n_t})^T$
	\item[$\pi^k$] A policy to be executed at time step $k$.
\end{IEEEdescription}

\section{Introduction}
\IEEEPARstart{B}{attery} energy storage systems (BESS) installed behind-the-meter have been increasingly popular at both residential and industrial consumers. This trend is mostly driven by decreasing prices, technological advancements and regulatory incentives. Increasing self-consumption from local generation 
by storing excess electricity generation for later use, is one of the major applications for installing behind-the-meter battery storage systems. 
For instance, in 2016, up to 46\% of installed PV systems smaller than 30~kWp were equipped with a battery storage system in Germany\cite{Speichermonitoring2017}.

In many cases however, the cost of a BESS remains high and the return on investment from solely self-consumption is too low\cite{Mulder2013}. Adding supplementary services to be delivered by the same BESS can lead to additional revenue streams and increase the return on investment.

A service for which BESS are deemed to be very appropriate is primary frequency control\cite{Knap2016} (also referred to as frequency containment reserve (FCR) or frequency response) due to their fast ramp rate\cite{Zhang2017}. Besides, frequency control has been identified as one of the highest value services for BESS \cite{Oudalov2006}. Complementary between primary frequency control and self-consumption can be expected, as primary frequency control is a service where power capacity is offered, while revenues from self-consumption are more driven by energy capacity. 

When participating in primary frequency control, one has to regulate his power consumption or injection to react proportionally on deviations of the grid frequency from the nominal frequency. The maximum contracted reserve capacity has to be activated when the frequency deviation is at a maximum, predefined value and within a predefined time limit. In the continental Europe (CE) synchronous area, this maximum is at a deviation of 200 mHz and has to reached within 30s\cite{ENTSO-E2013}.

In liberalized ancillary services markets, the TSO buys reserve capacity for frequency control from tertiary parties, who get paid for the power capacity they are able to sell.
In this work, we assume an end-consumer is able to offer this service to the TSO, possibly through an intermediary, for example an aggregator. 

When participating in primary frequency control, an asset has to be able to deliver the contracted power and follow the frequency signal during the whole contracted period. If the asset fails to do so, it is unavailable to provide the frequency control service and faces penalties charged by the TSO. As these penalties are usually high and the TSO expects an optimal service, in this work, we will try to constrain the risk of unavailability when delivering primary frequency control.

Being constantly available can be an issue when using energy-constrained assets such as a BESS. Over limited time periods, the frequency signal has a non-zero energy content and after having provided the service for a while, a BESS can be empty or full. In addition, efficiency losses in battery systems decrease their energy content or state of charge (SoC) when being charged and discharged continuously.
Therefore, a controller, or \emph{recharging strategy}, has to be in place to control the SoC to be within limits, ensuring that the reserve capacity remains available during the contracted period. 

Different strategies can be used: overdelivery (i.e. delivering more power than required), utilizing the deadband of the frequency signal (typically 10 or 20 mHz) to recharge or using a specific recharge controller that offsets the frequency control power for recharging.
A comparison of 
these methods is made in\cite{Fleer2016}, in which they conclude that overdelivery and deadband utilization alone is not sufficient to maintain the SoC within limits and an additional recharge controller is unavoidable. 
However, when using part of the BESS power for recharging, this part cannot be sold as a reserve capacity at the same time. One will thus have to optimize one's battery asset, maximizing the reserve capacity while minimizing the risk of unavailability.

The same is true when using the battery for the combination of frequency control and self-consumption. The BESS power used for self-consumption cannot be sold as reserve capacity for frequency control. One faces a trade-off between the two objectives which we incorporate into our optimization model. We see that complementarity between the two services occurs and show in a case study that optimally combining both services increases revenues from the BESS by 25\% compared to offering solely frequency control.

The remainder of the paper is organised as follows: in section~\ref{sec:rel_work}, related literature is reviewed and the main contributions of this paper are identified. The general problem treated in this paper is formulated in section~\ref{sec:model}. As it concerns a highly intractable problem, we treat the problem of providing solely frequency control first in section~\ref{sec:R1}. In section~\ref{sec:self-cons}, the objective of self-consumption is added, using a rule-based controller, optimized through stochastic optimization.
With the mathematical framework defined, section \ref{sec:simu} presents some simulations and results. Finally the paper is concluded and some suggestions for future work are given.


\section{Background and Related Work}\label{sec:rel_work}

\subsection{Background on Frequency Control}
While in Europe, liberalized markets exist for primary frequency control, this is not the case in North America. In North America, primary frequency control is traditionally delivered by generator governors or frequency responsive loads and is imposed as a requirement on large generators while no compensation is given for this service\cite{NERC-balancing}.

Markets do exist for regulation services, which is part of secondary frequency control, allowing third parties to offer their resources as regulation capacity.
Here, the asset has to follow a centrally dispatched signal to correct for the area control error (ACE) of the respective control area.
Compensation is not only based on offered capacity, but also on actual performance, rewarding assets that are able to perform better in following the regulation signal. 
Moreover, the California Independent System Operator (CAISO) has implemented a program for Non-Generator Resources (NGR) with Regulation Energy Management (REM) allowing for NGRs with limited energy content such as battery storage systems to competitively bid into the regulation market\cite{CAISO-REM}. PJM has implemented a high pass filter over its regulation signal in order to remove most of the energy content, making it more suitable for energy constrained resources such as BESS\cite{PJM}.

In the European context, some research has been conducted to create zero-mean frequency control signal\cite{Borsche2014}, but so far this has not been commercially implemented by any TSO.

While the approach presented in this paper can be applied to any type of frequency or regulation signal, the focus is on primary frequency control as defined by European TSOs, as they impose the strictest rules by requiring a 100\% availability and near perfect delivery. 

\subsection{Related Work}

From previous work on the provision of primary frequency control with a BESS, we identify two distinct approaches.
A first approach is to design a heuristic recharging strategy with simulations over historical frequency data for empirical optimization of the heuristic. 
For instance, Oudalev et al.\cite{Oudalov2007} design a rule-based recharge controller that acts when the SoC leaves the range $(SoC_{max}, SoC_{min})$. 
They use auxiliary resistors to consume additional power when the battery cannot provide enough, which we want to avoid in this work.
The heuristic recharging strategy presented in\cite{Borsche2013} is based on the moving average of the frequency signal, corrected for efficiency losses. The goal is to create a power profile with zero-mean, so that the battery does not get charged or discharged over time. 
A variant on this strategy is presented in\cite{Megel2013} and evaluated to give a higher return on investment when compared to the strategies from \cite{Oudalov2007,Borsche2013}. 
A rule-based control policy for fast energy storage unit in combination with a slower unit that is able to capture the energy content of the regulating signal is presented in\cite{Jin2011}.
While these heuristic strategies give good results, they do not ensure any form of optimality.

A second approach tries to overcome this by using more formal methods that can ensure optimality within the adopted framework. For instance, in \cite{Brivio2016}, a fuzzy control logic is used for primary frequency control and energy arbitrage in the Italian energy market. 
Zhang et al.\cite{Zhang2016} use dynamic programming to calculate an optimal recharging policy, recharging only when the frequency is in the deadband. 
Dynamic programming is also used for combining energy arbitrage and frequency regulation in the PJM regulatory zone \cite{Cheng2016}. However, both papers assume the reserve capacity a given parameter and are not able to optimize over this capacity itself. 


The combination of primary frequency control provision and minimization of photovoltaic (PV) and load curtailment by a battery storage system is considered in\cite{Megel2015}, where a model predictive control (MPC) is proposed to compute the allocation of the storage system for the two objectives. Although they model uncertainty in PV and consumption forecasts, they do not take this into account in the MPC controller.
Combination of self-consumption and primary frequency control is studied in  \cite{Braam2016}, however they use a heuristic controller that is only able to provide primary frequency control through pooling with a combined heat and power plant. Peak shaving and frequency regulation are combined in \cite{Shi2017}. By using the fast regulation signal from PJM, they are able to avoid the issue of limited energy content when offering frequency control services, which we want to overcome in this paper.

Using BESS connected to the distribution grid for frequency control might cause voltage problems or jeopardize the reliability of the distribution grid, when several BESS are connected to the same feeder\cite{Shahsavari}.
A potential solution can be a local voltage droop controller, which is shown to be effective in to avoid distribution grid constraint violations while having very limited impact on the performance of the service to be delivered\cite{DECONINCK2015}. 
In the remainder of this paper however, we will assume that the BESS providing frequency control are sufficiently dispersed over various feeders and do not endanger the reliability of the distribution grid.


Finally, it is worth mentioning that lately, there has been some commercial interest in the combination of self-consumption and frequency control with residential battery storage systems in Germany. More specifically, both companies \emph{Caterva}\cite{Caterva-RL} and \emph{Sonnen}\cite{Sonnen-RL} have presented a concept to combine self-consumption from PV with frequency control with a residential battery storage system. In both cases, the company acts as the intermediary party, operating a part of the storage systems for frequency control and offering the aggregated frequency control capacity to the TSO.

In this paper, we complement previous work by elaborating a controller that co-optimizes self-consumption, the reserve capacity and a recharge controller for primary frequency control. The main contributions can be summarized as follows:
\begin{itemize}
	\item We propose an optimized controller to maximize reserve capacity, which is able to provide more reserve capacity compared to the heuristic methods proposed in the literature.
	\item Building further on the work of Vrettos et al.\cite{Vrettos2016}, we extend their robust optimization approach towards a BESS application, and propose a new uncertainty set that provides explicit probability guarantees on battery constraint violation when providing frequency control.
	\item By co-optimizing self-consumption and frequency control, our approach is able to obtain more value than by using the BESS completely for only one of the objectives.
\end{itemize}

\section{Problem formulation}\label{sec:model}
We consider a simple, discrete battery model subject to a stochastic demand and production profile~$\bm{P}_{prof}~=~\bm{P}_{dem}~-~\bm{P}_{prod}$ and normalized frequency deviations $\bm{\Delta f}$. We model the BESS with constant charging and discharging efficiencies $\eta^c$, $\eta^d$. The battery has an effectively usable energy capacity in the range $(E_{min},E_{max})$ in which it is assumed to be able to provide the power range $(P_{min},P_{max})$.

The price for electricity injection into the grid is assumed smaller than the price for electricity consumption ${c}_{inj} < {c}_{cons}$, as this is imperative to make self-consumption financially interesting. 
The price for primary frequency control $c_r$ is assumed to be known, while the capacity $r$ is a variable to be optimized.

The objective of the problem is to minimize expected electricity consumption costs and maximize profits from primary frequency control, while keeping the BESS within its energy and power constraints. This results in following stochastic optimization program:
\begin{mini!}[1]
{}{\E[({c}_{cons} [\bm{P}_{grid}]^+ - {c}_{inj}[-\bm{P}_{grid}]^+) \Delta t] - c_r r } 
{\label{eq:extactprob}}{}
\addConstraint{\bm{P}_{grid}}{= \bm{P}_{prof} + \bm{P}_{bat}}
\addConstraint{\bm{P}_{bat}}{= \bm{P}_{ctrl} + r\bm{\Delta f}}
\addConstraint{E_{min}}{ \leq \bm{E}_{bat} \leq E_{max} \label{constr:Ebat}}
\addConstraint{P_{min} + r}{\leq \bm{P}_{ctrl} \leq P_{max} - r 
\label{constr:Pbat}}
\addConstraint{	E_{k+1}^{bat}}{ = E_{k}^{bat} + (\eta^{c} [P_{k}^{bat}]^+ - \frac{1}{\eta^{d}}  [-P_{k}^{bat}]^+) \Delta t. \label{constr:energy}}
\end{mini!}
Here, $\E[\cdot]$ denotes the expected value operator and $[\cdot]^+~\equiv~\max(\cdot,0),$ operating element-wise on vectors.
The power vector $\bm{P}_{grid}$ is the power that is actually put on the grid, consisting of the battery power $\bm{P}_{bat}$ and the demand profile $\bm{P}_{prof}$. 
The BESS power consist of two parts. One part is due to the frequency control and thus equal to the capacity times the frequency deviations $r\bm{\Delta f}$. A second part $\bm{P}_{ctrl}$ is dedicated to control the battery state of charge when providing frequency control, while optimizing the self-consumption.
Self-discharge losses, not incorporated here, can be added by subtracting them from the energy equation (\ref{constr:Ebat}).

To account correctly for the energy content of the battery~(\ref{constr:Ebat}), we assume that all power values are kept constant over one time step $\Delta t$. However, this is not possible when providing primary frequency control, as the BESS typically has to react within seconds to the real frequency control signal. Therefore, we define the discrete normalized frequency deviations $\Delta f_k$ as the average value over one time step:
\begin{equation*}
\Delta f_k = \frac{1}{\Delta t}\int_{(k-1)\Delta t}^{k\Delta t} \frac{(f(t)-f_{nom})}{\Delta f_{max}} dt,
\end{equation*} 
with $f(t)$ the real frequency, $f_{nom}$ the nominal value and $\Delta f_{max}$ the maximum frequency deviation on which on a has to react (for instance 200 mHz in the CE region).
To ensure that the instantaneous reserve capacity $r$ is always available, we have added it been explicitly to the hard power constraints in (\ref{constr:Pbat}).

As both the power profile $\bm{P}_{prof}$ and the frequency deviations $\bm{\Delta f}$ are stochastic variables that are gradually revealed over time, problem (\ref{eq:extactprob}) is a multi-stage stochastic program. This means that the ``here and now'' decision of the control power $\bm{P}_{ctrl}$ can be relaxed to a ``wait and see'' decision and depend on the past realisations of the power profile and frequency deviations ${P}_k^{ctrl} = \pi^k(P_1^{prof}, \ldots,P_k^{prof},\Delta f_1, \ldots,\Delta f_k)$ \cite{Shapiro2009}.
This is not true for the frequency control capacity $r$, as this value should be contracted with the TSO before the actual delivery takes place and one is not allowed to change this capacity during the delivery period.

Problem (\ref{eq:extactprob}) is a multi-stage non-linear stochastic program, which quickly becomes computationally intractable. To simplify, we propose to split the control power into two separate parts: a part for self consumption $ \bm{P}_{sc}$ and a part for recharging after frequency control activations $\bm{P}_{rc}$. Each is then depending on only one source of uncertainty:
\begin{equation*}
 \bm{P}_{ctrl} = \bm{P}_{sc}(\bm{P}_{prof}) + 
 \bm{P}_{rc}(\bm{\Delta f})
\end{equation*}

We can now look at (\ref{eq:extactprob}) as the combination of two distinct sub-problems: providing frequency control with a BESS and optimizing self-consumption. These sub-problems can then be put together, according to (\ref{eq:extactprob}), for joint optimization, which is expected to yield a better solution than the simple addition of the two objectives.


\section{Primary Frequency Control}\label{sec:R1}
In this section we will try to approximately solve problem~(\ref{eq:extactprob}), without the objective of self consumption (i.e.~$\bm{P}_{prof}~=~0$). 
The focus will be on the determination of the maximum frequency capacity $r$ the BESS can provided and the recharging policy $\bm{P}_{rc}(\bm{\Delta f})$ needed to keep the risk of unavailability as low as possible. 

\subsection{Recharging Policy}
The goal is to design a controller that ensures that the energy constraints (\ref{constr:Ebat}) are not violated when providing frequency control.
The typical control problem is to design a control policy which is a function of the current and past states of the system, here $P_{k}^{bat,rc} = \pi^k(E_1^{bat},\ldots,E_k^{bat})$. 
To come to a problem that can be solved efficiently, we will restrict ourselves to a linear policy. When writing this policy as an linear policy on the \emph{disturbance} $\bm{\Delta f}$ instead of the state, the problem becomes tractable~\cite{Ben-Tal2009}. The restriction to a disturbance feedback policy is not a limitation as it has been shown that a linear policy on the disturbance is as at least as general as an linear state feedback policy\cite{Goulart2006b}.
We can thus write the recharging policy as:
\begin{equation}\label{eq:recharge_policy}
P_k^{rc} = \sum_{i=1}^{k-1} d_{ki}\Delta f_i, \quad \bm{P}_{rc} = D\bm{\Delta f},  
\end{equation}
with $d_{ki}$ the coefficients of the recharge strategy, contained in the lower triangular matrix $D \in \mathbb{R}^{n_t \times n_t}$ with zeros on the diagonal. 
Note that we only sum up to $k-1$ in~(\ref{eq:recharge_policy}) so that there is no interference of the recharging power with the instantaneous frequency deviation $\Delta f(t)$.

One can interpret this recharging policy as a filter applied to the frequency control signal that creates a zero-mean signal, comparable to \cite{Megel2013,Borsche2013}. In this case, the recharge policy represents a change in the baseline on which the battery will perform the required frequency control activations. 

An aggregator can also pool the BESS together with another flexibility resource that is able to compensate for the recharging policy~\cite{RestorePatent}. Together they are able to follow the frequency signal exactly.

\subsection{Battery Efficiency}
Using the linear recharging policy (\ref{eq:recharge_policy}), problem~(\ref{eq:extactprob}) results in a mixed-integer stochastic program, which is known for its high computational complexity\cite{Dyer2006}. Therefore, we will use a heuristic approximation to turn~(\ref{eq:extactprob}) into a linear stochastic program. The integer variables in~(\ref{eq:extactprob}) arise purely because of the efficiencies $\eta^c, \eta^d$. By assuming an ideal battery and setting $\eta^c = \eta^d = 1$ in~(\ref{constr:energy}), the integer variables are removed and~(\ref{eq:extactprob}) becomes a linear problem.

As setting  $\eta^c = \eta^d = 1$ can be quite a coarse approximation, we instead apply the efficiencies to the frequency deviations:    
\begin{equation}\label{eq:en_freq}
\Delta f_k = \frac{1}{\Delta t}\int_{(k-1)\Delta t}^{k\Delta t} \Big(\eta^c\Big[\frac{\Delta f(t)}{\Delta f_{max}}\Big]^+ - \frac{1}{\eta^d}\Big[-\frac{\Delta f(t)}{\Delta f_{max}}\Big]^+\Big) dt,
\end{equation} 
which is exact if $\sign(\Delta f_k) = \sign(P_k^{bat,rc})$. By transforming the resulting disturbance feedback policy to an equivalent state feedback policy, it is possible to react appropriately to the impact of the efficiency.
Detailed simulations with real frequency data presented in section~\ref{sec:simu} demonstrate the validity of this approximation.

\subsection{Chance-Constraints and Robust Optimization}\label{sec:R1RO}

When applying the linear recharging policy from (\ref{eq:recharge_policy}), the power and energy content of the BESS are fully determined by the frequency deviations. The frequency deviation vector $\bm{\Delta f}$ is a multivariate stochastic variable in $\mathbb{R}^{n_t}$. This means that constraints~(\ref{constr:Ebat}), (\ref{constr:Pbat}) are actually probabilistic constraints, or so-called chance-constraints~\cite{Charnes1959}, 
and one has to constrain the probability of violation to be at maximum $\epsilon \in (0,1)$:
\begin{equation}\label{eq:chance_constr}
\Pr( \bm{a_i}^T \bm{\Delta f} \leq b_i ) \geq 1-\epsilon, \qquad i = 1,\ldots,n_c.
\end{equation}
Here, $n_c = 4n_t$ is the total number of constraints in~(\ref{constr:Ebat}), (\ref{constr:Pbat}) and $(\bm{a_i}, b_i)$ are defined to represent one constraint of~(\ref{constr:Ebat}),~(\ref{constr:Pbat}).

As breaching these constraints means that the frequency control service cannot be delivered, we want to make sure that the risk that this happens is as small as possible. Therefore, the goal is to get $\epsilon$ on the order of $10^{-4}$ or $10^{-5}$.

Several approaches to solve a chance-constrained problem exist. 
A popular approach is to use Monte Carlo sampling to approximate the real value of the probability in~(\ref{eq:chance_constr}). Explicit bounds on the number of samples are given in~\cite{Calafiore2006,Margellos2014} and are on the order of $O(n_{\delta}/\epsilon)$, with $n_{\delta}$ the dimension of the uncertainty. This would lead to a sample size on the order of $10^6$ for $\epsilon=10^{-4}$, which is not feasible if one considers a horizon of one day or more as we intend in this paper. Generating additional samples would require complete knowledge of the multivariate distribution of $\mathbf{\Delta f}$, which is never completely possible when working with observed data.

Analytic reformulation of (\ref{eq:chance_constr}) into a second-order cone constraint is possible if one assumes a Gaussian distribution\cite{Nemirovski2012}, which is not the case when considering $\bm{\Delta f}$.

One can also use a safe, convex approximation of (\ref{eq:chance_constr}). The conditional value-at-risk (CVaR)\cite{Rockafellar2000} is typically used as it is the tightest convex approximation to (\ref{eq:chance_constr}):

\begin{equation}\label{eq:CVAR}
\text{CVaR}_i^{1-\epsilon}(\bm{a_i}^T \bm{\Delta f} - b_i) \equiv \min_{\beta} 
\Big\{\beta + \frac{1}{\epsilon} 
\E  [\bm{a_i}^T \bm{\Delta f} - b_i -\beta ] ^+  \Big\} \leq 0,
\end{equation}
where $\E[\cdot]^+ \equiv \E[\max(\cdot,0)]$.
Despite its convexity, the CVaR risk measure is difficult to evaluate as the evaluation of $E[\cdot]^+$ requires multidimensional integration over the $\max(\cdot,0)$ operator. 
A sample average approximation of (\ref{eq:CVAR}) requires complete knowledge of the multivariate distribution and a large amount of samples to be accurate at small $\epsilon$ \cite{Shapiro2009}. 

Finally, the paradigm of robust optimization\cite{Ben-Tal2009} can be used to construct safe, tractable approximations to chance-constraints. The concept is to design an uncertainty set $\mathcal{F}$ of frequency deviations $\bm{\Delta f} \in \mathcal{F}$, against which the constraint has to be satisfied at all times:
\begin{equation*}
\bm{a_i}^T \bm{\Delta f} \leq b_i, \quad \forall \bm{\Delta f} \in 
\mathcal{F}, 
\qquad i = 1,\ldots,n_c.
\end{equation*}
This is equivalent to following worst-case formation:
\begin{equation}\label{eq:RO}
\max_{\bm{\Delta f} \in \mathcal{F}} \bm{a_i}^T \bm{\Delta f} \leq b_i, \qquad i = 1,\ldots,n_c.
\end{equation}

By correct design of $\mathcal{F}$, the solution of (\ref{eq:RO}) can ensure that the probability in (\ref{eq:chance_constr}) is bigger than or equal to the $(1-\epsilon)$ required, while retaining a tractable problem. 

Different uncertainty sets are proposed in literature 
(see e.g.\cite{Bertsimas2004,Ben-Tal2009,Chen2007}) 
of which their robust counterparts are shown to be upper bounds on the CVaR measure~(\ref{eq:CVAR})\cite{Chen2010}. Clearly, one is looking for the uncertainty set that provides the tightest upper bound. The robust counterparts of five types of tractable uncertainty sets are compared to the actual value of the CVaR measure in\cite{Chen2009} and it is shown that for small $\epsilon$, the uncertainty set based on forward and backward deviations provides the tightest bound to (\ref{eq:CVAR}). Furthermore, the robust counterpart is second-order cone representable and thus efficiently solvable by commercial solvers.

Unfortunately, the probability guarantee is only applicable to independent variables with zero mean, which is not the case when considering $\bm{\Delta f}$. However, by applying a whitening transformation\cite{kessy2017whitening}, one can obtain independent variables with zero mean:
\begin{equation}\label{whitening}
\bm{\tilde{f}} = W (\bm{\Delta f} - \overline{\bm{\Delta f}} ),
\end{equation}
where $\overline{\bm{\Delta f}}$ is the mean of $\bm{\Delta f}$, 
and $W^T W = \Sigma_{\bm{\Delta f}}^{-1}$ the Cholensky decomposition of the inverse of the covariance matrix $\Sigma_{\bm{\Delta f}}$ of ${\bm{\Delta f}}$.

The forward and backward deviations allow to include distributional asymmetry in the uncertainty set. They are defined for the stochastic variable 
$\tilde{f}_i$ as:
\begin{IEEEeqnarray*}{rCl}\IEEEnonumber
\sigma_{fi}(\tilde{f}_i) &=& \sup_{\theta>0} \sqrt{2\ln{\E[\exp(\theta \tilde{f}_i ) ] /  \theta^2}}, 				\IEEEnosubnumber \\
\sigma_{bi}(\tilde{f}_i) &=& \sup_{\theta>0} \sqrt{2\ln{\E[\exp(-\theta \tilde{f}_i ) ] / \theta^2}},   				\IEEEnonumber
\end{IEEEeqnarray*}
with $\theta \in \mathbb{R}$. The superior over $\theta$ can be found by applying a line search and approximating the expected value by its empirical average over the sampled data.
Define also $Q = \text{diag}(\sigma_{f1},\ldots,\sigma_{fn_t} )$ and $R = \text{diag}(\sigma_{b1},\ldots,\sigma_{bn_t})$. The uncertainty set $\mathcal{F_{\epsilon}}$ becomes then:
\begin{IEEEeqnarray*}{rCl}
\mathcal{F_{\epsilon}} &= \{&\mathbf{f}: \exists \bm{\beta}, \bm{\theta} \in \mathbb{R}_+^{n_t}, \mathbf{f} = \bm{\beta} -  \bm{\theta}, \nonumber \\ &&
 \lVert Q^{-1}\bm{\beta} + R^{-1}\bm{\theta} \rVert_2   \leq \sqrt{-2 \ln{\epsilon}}  \}
\end{IEEEeqnarray*}

Following \cite{Chen2010} and using (\ref{whitening}) to obtain independent variables $\bm{\tilde{f}}$ from $\bm{\Delta f}$, the $\text{CVaR}^{1-\epsilon}$ in (\ref{eq:CVAR}) is bounded by the worst-case of the constraint over the uncertainty set $\mathcal{F_{\epsilon}}$:
\begin{equation*}
\text{CVaR}_i^{1-\epsilon} \leq  \bm{a_i}^T \overline{\bm{\Delta f}} + \max_{\bm{\tilde f} \in \mathcal{F_{\epsilon}}} \bm{a_i}^T W^{-1} \bm{\tilde f} - b_i  \leq 0.
\end{equation*}
Finally, this can be reformulated as a second-order cone constraint\cite{Chen2009}:
\begin{equation}\label{eq:R1_constr}
\bm{a_i}^T \overline{\bm{\Delta f}} + \sqrt{-2\ln{\epsilon}} 
\lVert \bm{u_i}  \rVert_2
 \leq b_i, \qquad i = 1,\ldots,n_c ,
\end{equation}
where $\bm{u_i} = \max(Q\bm{a_i}^T W^{-1},-R\bm{a_i}^T W^{-1})$, with the maximum taken element-wise. Note that $\epsilon$ in (\ref{eq:R1_constr}) is under the logarithm, so that small values can easily be used.

\subsection{Equivalent State Feedback Policy}
The recharging strategy of (\ref{eq:recharge_policy}) is a disturbance feedback policy calculated with the efficiencies incorporated in the frequency signal (\ref{eq:en_freq}) and not in the battery model. This policy will therefore not be directly usable on a real battery system. However, by reformulating the policy as an equivalent state feedback policy it becomes practically usable. As efficiency losses are included in the state of charge of the battery, a state feedback policy can react on efficiency losses 
appropriately. 

Following\cite{Goulart2006b}, an equivalent state feedback policy can be calculated as:

\begin{equation}\label{eq:state_feedback}
\bm{P_{rc}} = (I+\frac{1}{r}D)^{-1} \frac{1}{r} D\bm{\Delta E_{bat}},
\end{equation}
with $\Delta E^{bat}_k = E_k^{bat}-E_{k-1}^{bat}$. In this form, the recharge power depends linearly on the past states, rather than on the past disturbances.

\section{Self-Consumption}\label{sec:self-cons}
In this section we will add the second part of the problem~(\ref{eq:extactprob}), i.e. finding a policy $\bm{P}_{sc}(\bm{P}_{prof})$ that allows to optimize self-consumption, while keeping in mind that a part of the battery has to be reserved for providing the primary frequency control. 

\subsection{Self-Consumption Policy}
The objective of self-consumption is to minimize the expected value of the total cost of electricity for the end-consumer. 
When facing constant consumption and production prices, a simple, rule-based control policy has proven to be very effective for this objective. The basic concept is to charge when there is more production than consumption and the battery is not full, and to discharge when there is more consumption than production and the battery is not empty.

To ensure sufficient energy and power of the BESS remains available for frequency control, we adapt the energy and power limits between which the battery can perform self-consumption to be smaller than the actual limits of the battery $(E_{max,k}^{sc} \leq E_{max}, E_{min,k}^{sc} \geq E_{min})$ and $(P_{max,k}^{sc} \leq P_{max}, P_{min,k}^{sc} \geq P_{min})$.
By making these limits dependent on the time $k$, they can be shaped towards the expected amount of generation or consumption. The control policy for self-consumption becomes then:
\begin{equation}
\label{eq:SC_rules}
\setlength{\nulldelimiterspace}{0pt}
P^{sc}_k=\left\{
	\begin{IEEEeqnarraybox}[\relax][c]{l?sc}
	\min (-P_k^{prof},P_{max,k}^{sc}),  & if & \left\{
		\begin{IEEEeqnarraybox}[\relax][c]{ll} 
			P_k^{prof} & < 0, \\ 
	 		E_k^{sc} & < E_{max,k}^{sc}, 
	 	\end{IEEEeqnarraybox} \right. \\
	\max (-P_k^{prof},P_{min,k}^{sc}), & if & \left\{
		\begin{IEEEeqnarraybox}[\relax][c]{ll} 
			P_k^{prof} & > 0, \\  
			E_k^{sc}  & > E_{min,k}^{sc} ,
		\end{IEEEeqnarraybox} \right. \\
	0,  & \multicolumn{2}{l}{otherwise.}
	\end{IEEEeqnarraybox}\right.
\end{equation}
This policy allows one to capture the most value from self-consumption while ensuring the capacity needed to deliver the frequency control is always available.

Notice that, when using this policy in combination with the frequency control policy described previously, one is actually dividing the battery into two virtual batteries with varying energy and power capacities: one for self-consumption and one for frequency control. 
Therefore, an estimation of the energy content of the virtual battery for self-consumption $E^{sc}$ should be available. This can be obtained by integrating $P^{sc}$, taking into account efficiency losses and other non-linearities as much as possible.
Alternatively, one can keep track of the energy content due to frequency control $E^{r}$ by integrating the corresponding power set-points $P^{rc}+r \Delta f$ and subtracting it from the measured state of charge: $E^{sc} = E^{bat} - E^{r}$.

\subsection{Stochastic Optimization}
Optimizing the self-consumption is a stochastic program in which the objective contains the expected value of the consumption and injection power vector:
\begin{mini}
{}{\E[({c}_{cons} \bm{P}_{cons} + {c}_{inj}\bm{P}_{inj}) \Delta t]}
{\label{eq:selfcons}}{}.
\end{mini}
A closed-form of this expected value is not readily available. Therefore, we will approximate the expected value by the sample average approximation (SAA)\cite{Shapiro2009}. By using various scenarios $j=1,\ldots, n_{sc}$ of the profile $\bm{P}_{prof}^j$, the empirical average of the objective will approximate the true expected value (\ref{eq:selfcons}).

By splitting the power for self-consumption into a part for charging and a part for discharging $\bm{P}_{sc} = \bm{P}_{sc,c} + \bm{P}_{sc,d}$, the efficiency can be accounted for correctly. As long as $c_{cons}>c_{inj}$, there is a cost for consuming energy and an optimal solution will always set $P_k^{sc,c} \cdot {P}_k^{sc,d} = 0, \forall k$.

Together with the constraints imposed by the rule-based charging policy in (\ref{eq:SC_rules}), one gets a linear program that can be solved efficiently:
\begin{mini}[1]
	{}{ \frac{1}{n_{sc}}\sum_{j=1}^{n_{sc}} ( {c}_{cons} \bm{P}_{cons}^j \Delta t + {c}_{inj} \bm{P}_{inj}^j\Delta t) }
	{\label{eq:selfcons_sampled}}{},
	\addConstraint{ \bm{P}_{cons}^j + \bm{P}_{inj}^j}{ = \bm{P}_{prof}^j + \bm{P}_{sc,c}^j + \bm{P}_{sc,d}^j}
	\addConstraint{0}{\leq \bm{P}_{sc,c}^j,\bm{P}_{cons}^j}
	\addConstraint{\bm{P}_{sc,d}^j,\bm{P}_{inj}^j}{\leq 0}
	\addConstraint{\bm{E}_{min}^{sc} \leq \bm{E}_{sc}^j}{\leq \bm{E}_{max}^{sc}}
	\addConstraint{\bm{P}_{sc,c}^j}{\leq \bm{P}_{max}^{sc}}
	\addConstraint{\bm{P}_{min}^{sc}}{\leq \bm{P}_{sc,d}^j}
	\addConstraint{E_{min} \leq \bm{E}_{min}^{sc}}{\leq \bm{E}_{max}^{sc} \leq E_{max}}
	\addConstraint{P_{min} \leq \bm{P}_{min}^{sc}}{\leq \bm{P}_{max}^{sc} \leq P_{max}}
	\addConstraint{E_{k+1,j}^{sc} }{=E_{k,j}^{sc} + (\eta^c {P}_{k,j}^{sc,c} + \frac{1}{\eta^d}P_{k,j}^{sc,d} ) \Delta t, }
\end{mini}
for all $j=1,\ldots,n_{sc}$ and $k = 1,\ldots, n_t$.
Here, we assume the scenarios or samples $\bm{P}_{prof}^j$ are independently identically distributed (iid). Samples with different probability distributions can be used by adding appropriate weights to each sample.

This problem can be combined with the chance-constrained problem of section~\ref{sec:R1} for providing frequency control, as shown in appendix \ref{app:A}, by adjusting the limits on energy content and BESS power in $b_i$ of~(\ref{eq:R1_constr}) to $(E_{max}-\bm{E}_{max}^{sc},\bm{E}_{min}^{sc}-E_{min})$ and $(P_{max}-\bm{P}_{max}^{sc},P_{min}-\bm{P}_{min}^{sc})$. 


\subsection{Scenario Reduction}\label{sec:sc_red}
Although the objective of the SAA problem (\ref{eq:selfcons_sampled}) converges to the true value~(\ref{eq:selfcons}) for $n_{sc} \rightarrow \infty$, the rate of convergence is on the order of $O_p (n_{sc}^{-1/2})$\cite{Shapiro2009}. A considerably large number of samples will thus be needed for sufficient accuracy. To limit the size of the problem and keep it tractable, scenario reduction techniques can be applied.  
We will use the backward scenario reduction of single scenarios of Dupa{\v{c}}ov{\'{a}} et al.\cite{Dupacova2003} based on the Kantorovich distance, since it has shown the best performance in our case.

\subsection{Evaluation of the Solution Quality}
As the objective of the SAA problem (\ref{eq:selfcons_sampled}) is an approximation to the true objective value~(\ref{eq:selfcons}),
 it would be instructive to have an estimation on how close the approximation is to the true value. 
Mak et al.\cite{Mak1999} provide a statistical method for calculating an upper and lower bound to the true objective value and the optimality gap of the SAA problem.

Let $\hat{x}$ be the optimal variables of the SAA problem, and $\bm{P}_{prof}^{j},\quad j=1,...,n_U$ iid profile samples, possibly different from the ones used in the SAA problem. 
Define $g(\hat{x},\bm{P}_{prof}^{j})$ as the objective of (\ref{eq:selfcons_sampled}) evaluated at $\hat{x}$ with $\bm{P}_{prof}^{j}$. 
An approximate $100(1-\alpha)\%$ confidence upper bound follows from the central limit theorem (CLT) on the average $\bar{g}_{n_U} = 1/n_U \sum^{n_U}_{j=1} g(\hat{x},\bm{P}_{prof}^{j})$ of the $n_U$ samples.

An $100(1-\alpha)\%$ confidence lower bound can be estimated by solving the SAA problem (\ref{eq:selfcons_sampled}) to optimality $n_L$ times: $\hat{g}_{n_{sc}}^1,\ldots, \hat{g}_{n_{sc}}^{n_L}$. The average  $\bar{g}_{n_L} = 1/n_L \sum^{n_L}_{i=1} \hat{g}_{n_{sc}}^i$ of the samples $\hat{g}_{n_{sc}}^i$ follows a $t$-distribution with $n_L-1$ degrees of freedom.
%
Finally, an $100(1-2\alpha)\%$ confidence upper bound to the optimality gap at $\hat{x}$ can be expressed as:
\begin{equation}\label{eq:opt_gap}
gap(\hat{x}) = \bar{g}_{n_U}  - \bar{g}_{n_L} + z_{\alpha}\frac{\hat{\sigma}_{n_U}}{\sqrt{n_U}} + t_{\alpha,n_L-1}\frac{\hat{\sigma}_{n_L}}{\sqrt{n_L}},
\end{equation}

with $z_{\alpha}=\Phi^{-1}(1-\alpha)$, where $\Phi(z)$ is the cumulative density function of the standard normal distribution. Here, $\hat{\sigma}_{n_U}$ is the sample standard deviation of the $n_U$ upper bound objective values $g(\hat{x},\bm{P}_{prof}^{j})$, $ t_{\alpha,n_L-1}$ the $\alpha$-critical value of the $t$-distribution with $n_L-1$ degrees of freedom and $\hat{\sigma}_{n_L}$ the standard deviation of the lower bound samples $\hat{g}_{n_{sc}}^i$.

\section{Simulation and Results}\label{sec:simu}

In this section we will present simulations and results of the mathematical program defined above. 
With the presented framework, we are able to draw some interesting conclusions about batteries providing frequency control and self-consumption. We will focus first on the robust optimization framework for frequency control only and then add the stochastic optimization for self-consumption.

In the simulations we consider a time horizon of one day, discretized in time steps of 15 minutes, so $n_t = 96$. A time horizon of one day is motivated by the daily cyclicality of consumption and PV production profiles. The time step of 15 minutes seems appropriate for the recharging policy when providing frequency control, as the regulations set by ENTSO\=/E state that a frequency deviation should be resolved within 15~minutes in the CE-region\cite{ENTSO-E2013}. The charging and discharging efficiencies are chosen to be $\eta^c=\eta^d=\sqrt{0.90}$, corresponding to a round trip efficiency of 90\%.
 
All optimizations are performed using the YALMIP\cite{Yalmip2004} toolbox with  Gurobi 7.0.2\cite{Gurobi} in MATLAB.
 
\subsection{Primary Frequency Control}

To assess the performance of primary frequency control, we use locally measured frequency data in the CE synchronous region with a resolution of 1 second over a period of three years (2014 - 2016). Missing data points are linearly interpolated up to 60~seconds. Days with remaining missing data points are removed from the data set, retaining 1091 complete days or samples. To test the performance of the approach on out-of-sample data points, we select 70\% of this data set randomly as training data, used to calculate $\sigma_{fi}$ and $\sigma_{bi}$, leaving 30\% for validation. 
The maximum frequency deviation on which the battery has to react $\Delta f_{max}$ in (\ref{eq:en_freq}) is set to $200~\text{mHz}$ as required in the CE synchronous region\cite{ENTSO-E2013}.

\subsubsection{Robust Optimization}
Consider a residential battery configuration of 10~kWh and 7~kW, with an initial charge of 5~kWh. Running the robust optimization as elaborated in section \ref{sec:R1RO}, with the chance of violating the battery constraints $\epsilon=10^{-4}$, gives a maximum reserve capacity of 6.37~kW. According to (\ref{constr:Pbat}) only 0.63~kW is to be used for recharging.
This reserve capacity is somewhat higher than what we found in literature, e.g. \cite{Borsche2013} gives a maximum reserve capacity of 4.66~kW for the same battery configuration.

\begin{figure}[!t]
\centering
\includegraphics[width=\columnwidth]{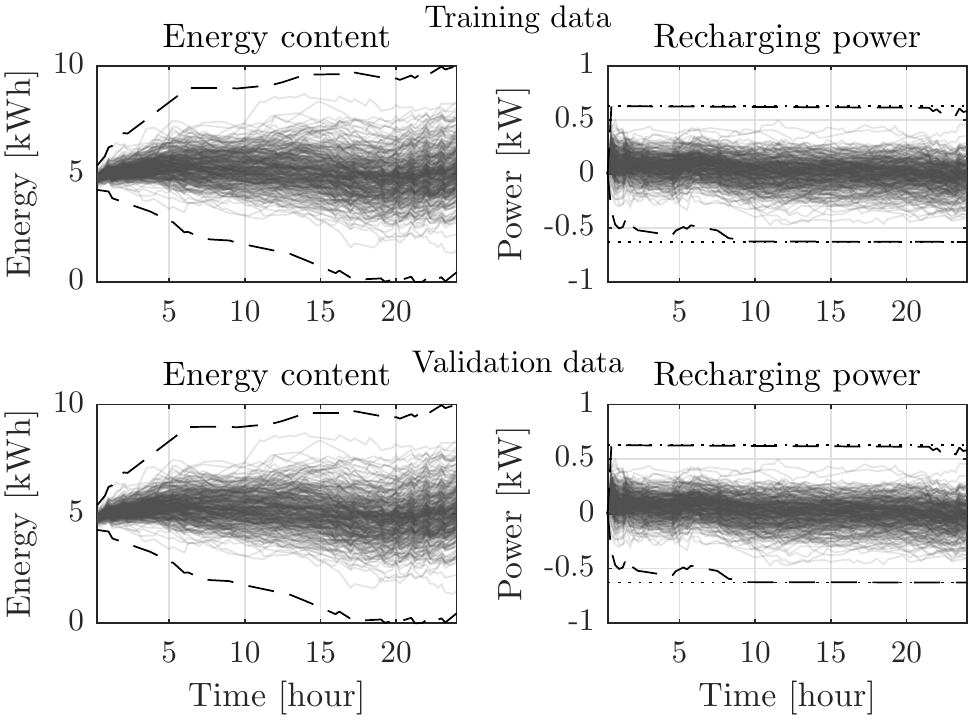}
\caption{Energy content and recharging power of the battery, for the frequency data in the training set (top) and the validation set (bottom), discretized according to~(\ref{eq:en_freq}) The dashed lines show the maximum and minimum cases, according to~(\ref{eq:RO}). The dotted lines show the maximum and minimum recharging power that is allowed, following (\ref{constr:Pbat}). Each line represents a frequency sample of one day.}
\label{fig:r1only}
\end{figure}

Figure \ref{fig:r1only} shows the corresponding energy and recharging power profiles for each frequency profile in the dataset. One can see that for both the training and the validation data, the energy content and recharging power stay well below the limits. Having chosen a small value of $\epsilon$ this makes sense, as the battery should be able to withstand more extreme frequency profiles that are not presented in the data sets. 
Using the uncertainty set $\mathcal{F_{\epsilon}}$ in~(\ref{eq:RO}), it is possible to calculate the maximum and minimum energy content and recharging power of the battery, shown by the dashed lines. One can see that they do not breach but do approach the boundary conditions of the battery, as expected.

\subsubsection{Equivalent State Feedback}

\begin{figure}[!t]
	\centering
	\includegraphics[width=0.9\columnwidth]{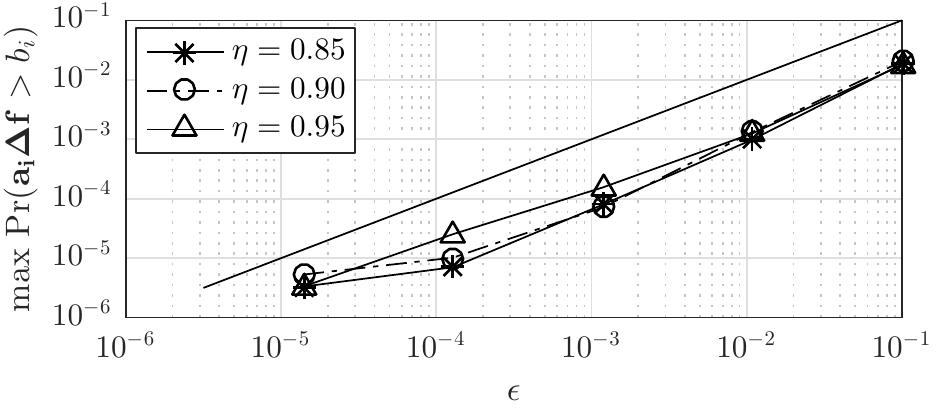}
	\caption{Maximum probability of constraint violation with the equivalent state feedback controller (\ref{eq:state_feedback}), for various values of $\epsilon$ and round trip efficiencies $\eta$, calculated with $n_R=10^6$ samples. The solid black line denotes the maximum allowed probability.}
	\label{fig:r1statefeedback}
\end{figure}
The results presented in figure~\ref{fig:r1only} use an ideal battery model without losses but with the efficiencies $\eta$ included into the frequency disturbances, as in (\ref{eq:en_freq}). 
To evaluate the performance of the state feedback controller of~(\ref{eq:state_feedback}) we have calculated the maximum probability of constraint violation:
\begin{equation}\label{eq:constraint_violation}
\text{max}_i \Pr (\mathbf{a_i \Delta f} > b_i)
\end{equation}
with the state feedback controller on a battery with a round-trip efficiency $\eta < 1$ for various values of $\epsilon$ and $\eta$ as follows.

To obtain the averaged frequency signal $\bm{\Delta f_k^o}$ that is not corrected for efficiency losses, we use (\ref{eq:en_freq}) with $\eta^c=\eta^d=1$. By then applying the whitening transformation (\ref{whitening}) on $\bm{\Delta f_k^o}$, we obtain independent variables with zero mean $\tilde{f}^o_k$, from which we can generate new frequency samples $\bm{\Delta{f}^r}$ by resampling ${\tilde{f}^o}_k$ with replacement $n_R$ times and applying the inverse of the whitening transformation. Using the state feedback controller~(\ref{eq:state_feedback}) with $\bm{\Delta{f}^r}$ gives then a Monte Carlo estimate of (\ref{eq:constraint_violation}) with $n_R$ samples.

Figure~\ref{fig:r1statefeedback} shows the resulting 99\% confidence upper bound (calculated according to p.217 in\cite{Shapiro2009}) of the maximum probability of constraint violation (\ref{eq:constraint_violation}) for various values of $\epsilon$ and $\eta$ with $n_R=10^6$ Monte Carlo samples. One can see that the actual probabilities stay well below the maximum allowed $\epsilon$, for all evaluated values of $\epsilon$ while the effect of $\eta$ is minimal.
%
	

\subsubsection{Maximum Reserve Capacity}
\begin{figure}[!t]
	\centering
	\includegraphics[width=0.9\columnwidth]{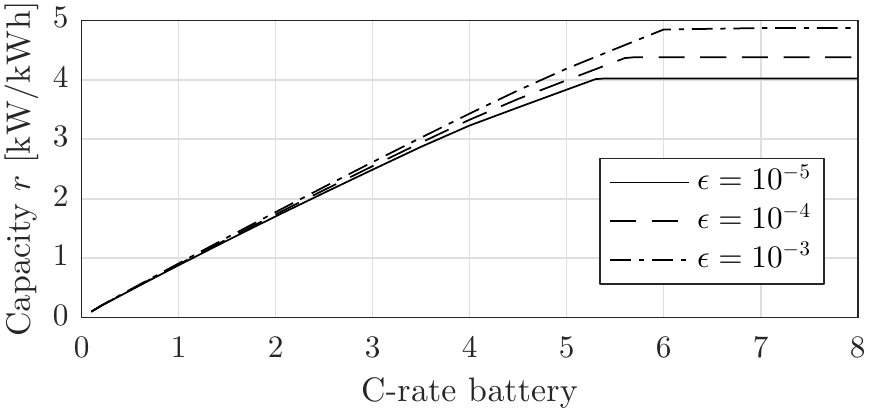}
	\caption{Primary frequency control capacity $r$ per kWh of battery capacity in function of the C-rate of the battery, for various values of $\epsilon$ and a round trip efficiency $\eta = 90\%$.}
	\label{fig:R1caps}
\end{figure}

When offering primary frequency control, it is the reserve capacity $r$ that is creating value for the BESS. Using the BESS solely for frequency control, the reserve capacity depends on the \emph{C-rate} of the BESS, defined as the maximum power divided by the maximum energy content, and the energy capacity of the BESS.

Figure~\ref{fig:R1caps} shows the maximum reserve capacity $r$ per kWh in function of the C-rate, for various values of $\epsilon$. As one can see, the relative reserve capacity is a concave function of the C-rate. Increasing the C-rate of a battery while keeping the energy content constant will thus increase the reserve capacity one can offer with this battery. This is an interesting result, as the main cost driver for batteries is the energy content, rather than the maximum power capacity.

The reserve capacity increases with the C-rate up to a maximum point, at which it is limited solely by the energy content of the battery. Increasing the maximum power of the battery beyond this point will not have any effect on the reserve capacity one can offer. The recharge policy is at its maximum, immediately compensating for the effect of the past frequency deviation. An increase in battery power will not have an effect any more on the recharge policy, thus not be able to increase the reserve capacity.

As could be expected, increasing the probability of battery constraint violation $\epsilon$ also increases the amount of reserve capacity one can offer with the same battery. However, this also means an increased risk of unavailability and penalties. 
If the battery is part of a pool of an aggregator, a higher $\epsilon$ can be chosen if the pool can be used as back-up when the BESS constraints are reached.

\subsection{Combination with Self-Consumption}
To asses the performance for the combination of frequency control and self-consumption we consider the same battery configuration as before. 
Residential demand profiles are generated from the CREST demand model\cite{Richardson2010} for a weekday in March. PV profiles are generated from the model presented in\cite{Bright2015} and scaled to represent a PV system of 4.0~kWp. We assume $c_r = 14.71$~EUR/MW/h, which was the average price for primary frequency control on Regelleistung in 2016\cite{Regelleistung}, $c_{cons} = 28.73$~cEUR/kWh, corresponding to the average consumption price in Germany in 2016\cite{BDEW_Strompreis}, and $c_{inj} = 12.20$~cEUR/kWh, the current Germany feed-in tariff for residential PV\cite{EEG2017}.

\subsubsection{Selection of Number of Scenarios}
%

Calculating the optimality gap using (\ref{eq:opt_gap}) with $n_U = 10^5, n_L = 10$ and $\alpha = 0.005$ for various numbers of scenarios $n_{sc}$, we find that overall, the optimality gap decreases rapidly to a small value ($\leq 3\%$ if $n_{sc}\geq250$) and as from about 1000 scenarios, the optimality gap can be expected to be less than 1\%.

When using the scenario reduction method from section~\ref{sec:sc_red}, an optimality gap smaller than 1\% can be reached from about 500 reduced scenarios.

\subsubsection{Self-Consumption and Frequency Control}
\begin{figure}[!t]
	\centering
	\includegraphics[width=\columnwidth]{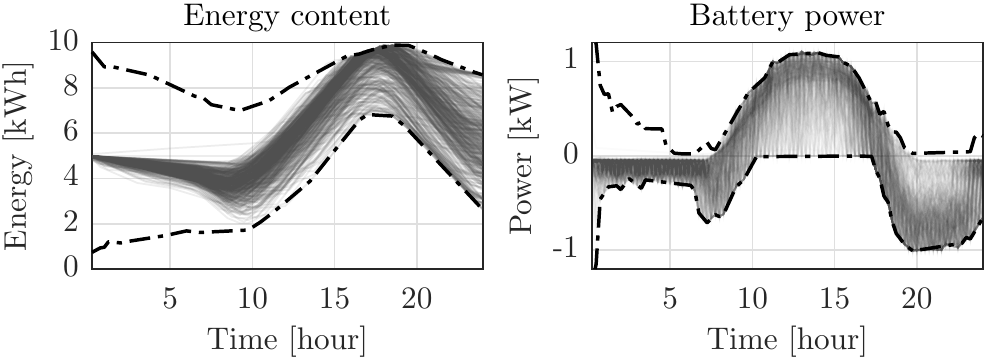}
	\caption{Energy content and battery power for self-consumption of 500 scenarios. The dashed black lines denote the limits $(\bm{E}_{min}^{sc},\bm{E}_{max}^{sc})$ and $(\bm{P}_{min}^{sc},\bm{P}_{max}^{sc})$ from the rule-based charging policy.}
	\label{fig:SC_R1}
\end{figure}
Combining primary frequency control with self-consumption with 500 scenarios from the scenario reduction gives an optimal reserve capacity of 5.65~kW. The remaining power (1.35~kW) is used for maximizing self-consumption and recharging for frequency control.
Figure \ref{fig:SC_R1} shows the BESS energy and power profiles for self-consumption of 500 scenarios. The dashed black lines denote the limits $(\bm{E}_{min}^{sc}, \bm{E}_{max}^{sc})$ and $(\bm{P}_{min}^{sc}, \bm{P}_{max}^{sc})$ from the rule-based charging policy (\ref{eq:SC_rules}).

At moments when production is expected to be high, during noon, the controller reserves power and energy in the battery to charge for the self-consumption objective, which can be discharged at times when expected consumption is higher, mainly in the evening. Less power is reserved during the night, as less consumption is expected at these times.

The expected value of self-consumption during this day is 0.81~EUR, while from frequency control with $r=5.65$~kW capacity at 14.71~EUR/MW/h, revenues are 2.00~EUR. In total, this gives a value of 2.81~EUR. When using the BESS only for self-consumption, the expected value is only 0.94~EUR. When using the BESS only for frequency control, the reserve capacity is a bit higher: $r=6.37$~kW, and total revenues are 2.25~EUR. The revenues of the combined optimization are thus more about 3 times higher compared to the case of only self-consumption and 25\% higher compared to solely frequency control.

\begin{figure}[!t]
	\centering
	\includegraphics[width=0.9\columnwidth]{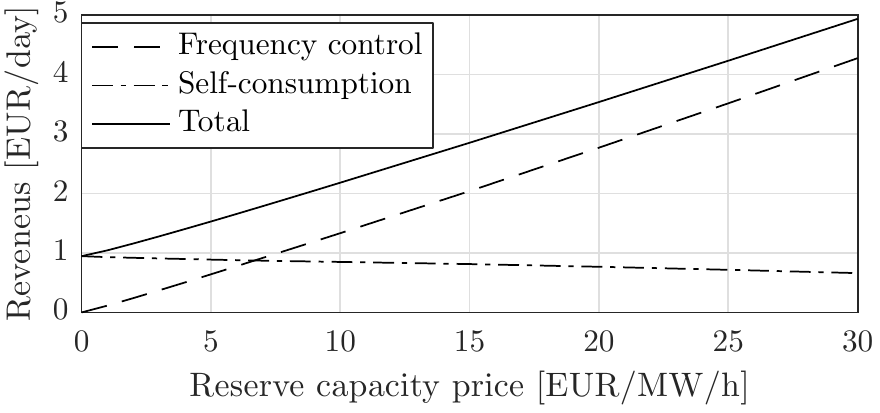}
	\caption{Revenues from combined frequency control and self-consumption in function of the reserve capacity price for the simulated day in March.}
	\label{fig:rev_R1prices}
\end{figure}

Figure~\ref{fig:rev_R1prices} shows the total revenues and the break-down into a part from self-consumption and a part from frequency control, for varying R1 prices. A trade-off between frequency control and self-consumption is clearly visible, as with increasing reserve capacity prices a larger part of the battery is reserved for frequency control and consequently, revenues from self-consumption decrease.

\section{Conclusion}\label{sec:conc}
In this paper, we have presented a framework for stochastic co-optimization of primary frequency control and self-consumption with a battery energy storage system.

Using robust optimization as a safe, tractable approximation to chance-constraints, we have design a linear recharging strategy that allows one to make the risk of unavailability arbitrarily small, while keeping the problem efficiently solvable. Simulations with real frequency data show the performance of the recharging controller when converted to an equivalent state feedback controller.

We have adopted a rule-based controlled to maximizing self-consumption, which allows to reserve more energy and battery power on moments when expected consumption or production is higher, while using other moments for recharging from primary frequency control. 
A sample average approximation is used to estimate the expected value of self-consumption and perform the trade-off between self-consumption and primary frequency control.

We have performed a case study on a residential battery system. The results show that there is a clear complementary in combining frequency control with self-consumption. 
Our co-optimization adds 25\% of value compared to the use for frequency control alone, while increasing value times 3 when compared to self-consumption alone.

Interesting future work is to look into non-linear recharging policies for frequency control, e.g. by direct policy search. 
Incorporating a more accurate battery model, where efficiencies and power limits can be dependent on the state of charge will represent reality better. One could also incorporate the costs of energy when providing frequency control into the model in a stochastic manner, in analogy to the way we treated self-consumption. 
While this paper only considers self-consumption and frequency control, other services, such as peak shaving, time of use tariff optimization or voltage control could also be incorporate into the optimization model.
Finally, validation of the battery control strategies on a real battery system should be performed.



\appendices
\section{Combined Optimization Problem for Frequency Control and Self-Consumption}\label{app:A}
The complete second-order cone program, combining frequency control and self-consumption (\ref{eq:selfcons_sampled}) is presented below:

\begin{mini*}
{}{ \frac{1}{n_{sc}}\sum_{j=1}^{n_{sc}} ( {c}_{cons} \bm{P}_{cons}^j \Delta t + {c}_{inj} \bm{P}_{inj}^j\Delta t) - c_r r }
{\label{eq:full_prob}}{},
\addConstraint{ \bm{P}_{cons}^j + \bm{P}_{inj}^j}{ = \bm{P}_{prof}^j + \bm{P}_{sc,c}^j + \bm{P}_{sc,d}^j}
\addConstraint{0}{\leq \bm{P}_{sc,c}^j,\bm{P}_{cons}^j}
\addConstraint{\bm{P}_{sc,d}^j,\bm{P}_{inj}^j}{\leq 0}
\addConstraint{\bm{E}_{min}^{sc} \leq \bm{E}_{sc}^j}{\leq \bm{E}_{max}^{sc}}
\addConstraint{\bm{P}_{sc,c}^j}{\leq \bm{P}_{max}^{sc}}
\addConstraint{\bm{P}_{min}^{sc}}{\leq \bm{P}_{sc,d}^j}
\addConstraint{E_{min} \leq \bm{E}_{min}^{sc}}{\leq \bm{E}_{max}^{sc} \leq E_{max}}
\addConstraint{P_{min} \leq \bm{P}_{min}^{sc}}{\leq \bm{P}_{max}^{sc} \leq P_{max}}
\addConstraint{E_{k+1,j}^{sc} }{=E_{k,j}^{sc} + (\eta^c {P}_{k,j}^{sc,c} + \frac{1}{\eta^d}P_{k,j}^{sc,d} ) \Delta t }
\addConstraint{\sqrt{-2\ln{\epsilon}} 
	\lVert \bm{u_i} \rVert_2}{\leq b_i -\bm{a_i}^T \overline{\bm{\Delta f}} }
\addConstraint{Q\bm{a_i}^T W^{-1}}{\leq \bm{u_i}}
\addConstraint{-R\bm{a_i}^T W^{-1}}{\leq \bm{u_i},}
\end{mini*}
for all $j=1,\ldots,n_{sc}$, $k = 1,\ldots, n_t$ and $i = 1,\ldots n_c$. If we define constraint matrix 
 $A = [D^T | -D^T  (D+rI)^TG^T | -(D+rI)^TG^T]^T$, with $G$ a lower triangular matrix with $\Delta t$ as elements, and vector $\bm{b} = [ P_{max}-{\bm{P}_{max}^{sc}}^T-r | P_{min}-{\bm{P}_{min}^{sc}}^T+r | E_{max}-{\bm{E}_{max}^{sc}}^T | {\bm{E}_{min}^{sc}}^T-E_{min}] ^T$, then $\bm{a_i}^T$ is the $i$-th row of $A$ and $b_i$ the $i$-th element of $\bm{b}$.
\ifCLASSOPTIONcaptionsoff
  \newpage
\fi



\bibliographystyle{IEEEtran}
\bibliography{IEEEabrv,bibl}{}

\end{document}